\documentclass[10pt]{amsart}
\usepackage{palatino}
\usepackage{amsfonts}
\usepackage{amssymb}
\usepackage{amscd}
\usepackage{hyperref}
\usepackage{graphicx}
\usepackage{tikz}

\newtheorem{thm}{Theorem}[section]

\newtheorem{lemma}[thm]{Lemma}
\newtheorem{prop}[thm]{Proposition}
\theoremstyle{definition}
\newtheorem{defn}[thm]{Definition}

\theoremstyle{rem}
\newtheorem{rem}[thm]{Remark}
\numberwithin{equation}{subsection}

\renewcommand{\emptyset}{\varnothing}

\newcommand{\etc}{\,\ldots}

\begin{document}

\title[The Chow ring of the moduli space of curves of genus zero]{The Chow ring of the moduli space of curves of genus zero}
\author[M. Tavakol]{Mehdi Tavakol}
   \address{Department of Mathematics,
              KTH, 100 44 Stockholm, Sweden}
   \email{tavakol@math.kth.se}

\maketitle
\centerline {\bf Introduction.}
\vskip 1cm

\subsection{Overview}
The purpose of this note is to give a presentation of the intersection ring of the moduli space $\overline{M}_{0,n}$ of stable $n$-pointed curves of genus zero. In the first section we recall general facts about the moduli space $\overline{M}_{0,n}$. There are several constructions of this space. We recall three of them and present another construction of $\overline{M}_{0,n}$ as a blow-up of the variety $(\mathbb{P}^1)^{n-3}$. In the next section we describe the intersection ring of the moduli space in terms of generators and relations. 
We give a description of a basis for the Chow groups. We show that with respect to our basis for the Chow groups there is an explicit  duality between the Chow groups in complementary degrees. 
This duality was implied by the result of Keel [K] from his inductive argument. He shows that dual bases for $A^*(\overline{M}_{0,n})$ induce dual bases for $A^*(\overline{M}_{0,n+1})$. This gives a \emph{recursive} way to get dual bases. Our presentation of the intersection 
ring is simpler in the sense that there are fewer generators and fewer relations, and we give an \emph{explicit} basis for the Chow groups. 

\subsection{Proof steps}
While there are well-known formulas for the intersection ring of blow-ups we use the formula of Keel \eqref{keel}, which simplifies the 
computations in our situation. This formula determines the intersection ring of certain blow-ups in terms of generators and simple relations.
We describe the space of relations in $A^*(\overline{M}_{0,n})$ in \ref{relations}. This leads to a basis for the Chow groups and the notion of \emph{standard monomials} introduced in Section \ref{standard}. In Proposition \ref{st} we see that Chow groups are  \emph{additively} generated by standard monomials. Then we observe that there is an involution on the Chow ring which gives a one-to-one
correspondence between standard monomials in complementary degrees. The matrices of the intersection pairing between standard monomials and their duals become of a simple form with respect to a natural ordering of standard monomials. 
More precisely, we show that the intersection matrices 
consist of identity matrices along the main diagonal up to a sign and all blocks below the main diagonal are zero.
The filtration \ref{filter} of the Chow ring together with the property \ref{fil} gives a geometric justification of the vanishing result formulated in Proposition \ref{tri}. To prove our claim about the square blocks on the main diagonal we use the identity \eqref{Number} in the Chow ring of $\overline{M}_{0,n}$. This identity is easily proven using the intersection theory of blow-ups in our construction of the moduli space. This completes the description of the Chow ring and the intersection pairings as it is stated in Theorem \ref{main}.
The last section contains examples. 

\vspace{+10pt}
\noindent{\bf Acknowledgments.}
The author thanks Carel Faber for useful discussions. 
\section{Construction of the moduli space $\overline{M}_{0,n}$}

Recall that the moduli space $M_{0,n}$ parametrizes the isomorphism classes of pointed curves $(C;x_1,\etc,x_n)$, where $C$ is a smooth curve of genus zero and the $x_i$'s are distinct points on $C$. Every smooth curve of genus zero
is isomorphic to the projective line $\mathbb{P}^1$ and for every 3 distinct points $a,b,c$ on $\mathbb{P}^1$ there exists a unique automorphism of the projective line sending these points to $0,1,\infty$ respectively. This means that the space $M_{0,n}$ coincides with $$U:=(\mathbb{P}^1-\{0,1,\infty\})^{n-3}-\Delta,$$
where $\Delta$ stands for the big diagonal consisting of all points where at least two coordinates become equal to each other. 
The moduli space of stable $n$-pointed curves of genus zero is introduced by 
Grothendieck in [D] as a compactification of the moduli space $M_{0,n}$ and is studied by Knudsen in [Kn].
This coincides with the Deligne-Mumford compactification of the moduli space
$M_{0,n}$ which is defined in a more general setting. We recall the definition.
Let $S$ be a scheme, and let $n \geq 3$ be an integer. 

\begin{defn}
A stable $n$-pointed curve of genus zero over $S$ is a flat and proper morphism $\pi:C \rightarrow S$ together with $n$ distinct sections $\sigma_i:S \rightarrow C$ such that

\begin{itemize} 
\item The geometric fibers $C_s$ of $\pi$ are reduced and connected curves with at most ordinary double points.

\item $C_s$ is smooth at $P_i=\sigma_i(s)$ for $1 \leq i \leq n$.

\item The points $P_i$ and $P_j$ are distinct for $i \neq j$.

\item The number of points where a nonsingular component $E$ of $C_s$ meets the rest of $C_s$ plus the number of points $P_i$ on $E$ is at least 3.

\item The curve $C_s$ is of arithmetic genus zero, i.e. $H^1(C_s,\mathcal{O}_{C_s})=0$. 
\end{itemize}
\end{defn}

Let $\overline{M}_{0,n}$ be the functor which sends a scheme $S$ to the collection of stable $n$-pointed curves of genus zero over $S$ modulo isomorphisms. 
In [Kn] Knudsen shows that this functor is represented by a smooth complete variety $X_n$ together with a universal curve $\pi: U_n \rightarrow X_n$, and universal sections $\sigma_1,\etc,\sigma_n$.  
His method is inductive: Suppose $U_n \rightarrow X_n$ with sections $\sigma_1,\etc,\sigma_n$ represents $\overline{M}_{0,n}$. Then $$U_n \times_{X_n} U_n \rightarrow U_n$$ with the pulled back sections 
$\sigma_1,\etc,\sigma_n$ and the additional section coming from the diagonal is the universal $n$-pointed curve with an additional section. There is a blow-up $(U_n \times_{X_n} U_n)^s$ of 
$U_n \times_{X_n} U_n$ (which he calls the stabilization of the family $U_n \times_{X_n} U_n \rightarrow U_n$) so that $$(U_n \times_{X_n} U_n)^s \rightarrow U_n$$ is the universal $(n+1)$-pointed curve. In particular, $X_{n+1}=U_n$.

As Keel remarks in [K] the computation of the intersection ring of $\overline{M}_{0,n}$ from this construction becomes difficult since the fibered product $U_n \times_{X_n} U_n$ is not smooth, and the blow-up is not 
along a regularly embedded subscheme. 

In [K] Keel gives an alternative inductive method. He shows that $\pi$ can be factored as 

$$\begin{CD}
X_{n+1} @>\rho>> X_n \times \mathbb{P}^1 \\ 
@. @VV\pi_1 V \\ 
@. X_n
\end{CD}$$
where $\pi_1$ is projection on the first factor, and $\rho$ is a composition of blow-ups of smooth varieties along smooth codimension two subvarieties. This enables him to give a complete description of the intersection ring
of the moduli space. In particular, he proves that the canonical map from the Chow groups to homology (in characteristic zero) is an isomorphism. He gives a recursive formula for the Betti numbers of $\overline{M}_{0,n}$. 
Keel also gives an inductive recipe for determining dual bases in the Chow rings and shows that the Chow ring is generated by divisors. He also describes the ideal of relations. 

M. M. Kapranov [Ka] gives another construction of $\overline{M}_{0,n}$ via a sequence of blow-ups:
Choose $n-1$ generic points $q_1,\etc,q_{n-1}$ in $\mathbb{P}^{n-3}$. The variety $\overline{M}_{0,n}$ can be obtained from $\mathbb{P}^{n-3}$ by a series of blow-ups of all the projective spaces spanned by the $q_i$.
The order of these blow-ups can be taken as follows:
\begin{enumerate}
\item[(1)] The points $q_1,\etc,q_{n-2}$ and all the projective subspaces spanned by them in order of increasing dimension.
\item[(2)] The point $q_{n-1}$, all the lines $\langle q_1,q_{n-1} \rangle,\etc, \langle q_{n-3},q_{n-1} \rangle$ and the subspaces spanned by them in order of increasing dimension.
\item[(3)] The line $\langle q_{n-2},q_{n-1} \rangle$, the planes $\langle q_i,q_{n-2},q_{n-1} \rangle, i \neq n-3$ and all subspaces spanned by them in order of increasing dimension, etc.
\end{enumerate}

We give another construction of $\overline{M}_{0,n}$ as a natural compactification of the moduli space $M_{0,n}$: Let $X$ be the $(n-3)$-fold product of the projective line $\mathbb{P}^1$ and $I$ 
be a subset of the set $\{1,\etc,n\}$ such that $$(*) \qquad |I \cap \{n-2,n-1,n\}| \leq 1 \qquad  \mathrm{and} \ 3 \leq |I| \leq n-2.$$ The subvariety $X_I$ of $X$ is defined as follows: 
If the intersection $I \cap \{n-2,n-1,n\}$ is the empty set, it is defined to be the set of all points $(x_1,\etc,x_{n-3})$ in $X$ where the coordinates corresponding to the index set $I$ are equal to each other. 
Otherwise, consider the following correspondence $$n-2 \leftrightarrow 0,  \qquad n-1 \leftrightarrow 1,  \qquad n \leftrightarrow \infty.$$
We refer to $n-2,n-1,n$ as \emph{special elements}. In this case the subvariety $X_I$ consists of all points in $X$ where the coordinates corresponding to non-special elements in $I$ are equal to the point specified via the above correspondence
to the special element of $I$. Notice that the subvariety $X_I$ is of dimension $n-2-|I|$.

For example, let $n=6$, $I=\{1,2,3\}$ and $J=\{2,3,4\}$. Then $$X_I=\{(x,x,x): x \in \mathbb{P}^1\}, \qquad X_J=\{(x,0,0): x \in \mathbb{P}^1\}.$$

We identify $M_{0,n}$ with the open subset $U$ of $X$ mentioned above consisting of all points whose coordinates are distinct from each other and don't belong to the set $\{0,1,\infty\}$. 
The universal family over $M_{0,n}$ is the trivial family $$\pi: M_{0,n} \times \mathbb{P}^1 \rightarrow M_{0,n},$$ where $\pi$ is the projection onto the first factor. Let the universal sections of this family be $\sigma_1,\etc,\sigma_n$. 
We assume that the last three sections are constant, sending a point $P$ on $M_{0,n}$ to the points $0,1,\infty$ of the fiber over $P$ respectively. 

The space $\overline{M}_{0,n}$ is obtained from $X$ by the following sequence of blow-ups:  First blow-up three points $X_{\{1,\etc,n-3,n-2\}}, X_{\{1,\etc,n-3,n-1\}}$ and $X_{\{1,\etc,n-3,n\}}$. This separates all lines 
$X_I$ for which $|I|=n-3$. Then blow up the proper transform of all these lines. We continue this process and increase the dimension of the blow-up centers at each step. Note that at each step the blow-up center is the regularly embedded union of the proper transform of subvarieties $X_I$. This process goes on as long as $|I| \geq 3$. The exceptional divisor of the blow-up along the 
subvariety $X_I$ is denoted by $D_I$ as well as the class of its proper transform under later blow-ups.

The family $\pi$ extends to a family over the resulting space $\widetilde{X}$ and permits $n$ disjoint sections in the smooth locus of the fibers. This requires a blow-up of the total space of the family along the inverse image via $\pi$ of the ideal sheaf of the blow-up centers on the base of the family. For the precise statement see Corollary 7.15 in [H]. We denote the resulting family by the same letter $\pi$ and its sections by $\sigma_1,\etc,\sigma_n$, by abuse of notation. The family $\pi$ induces a morphism $$F:\widetilde{X} \rightarrow \overline{M}_{0,n}.$$
The morphism $F$ sends the point $P \in \widetilde{X}$ to the isomorphism class of the pointed curve $$(\pi^{-1}(P);\sigma_1(P),\etc,\sigma_n(P)).$$ The restriction of the morphism $F$ on the open subset $U \subset \widetilde{X}$ defines an isomorphism between $U$ and the image $F(U)=M_{0,n}$. This means that $F$ is a birational morphism. The study of the fibers of $F$ and Theorem 7.17 in [H] show that $F$ is an isomorphism. 

\begin{rem}
There are standard divisor classes $D_I$ on the Deligne-Mumford compactification $\overline{M}_{0,n}$ of $M_{0,n}$ whose union is the boundary $\overline{M}_{0,n} -M_{0,n}$. These boundary divisors are indexed by subsets 
$I \subset \{1,\etc,n\}$, where $2 \leq |I| \leq n-2$. The generic curve parameterized by the divisor $D_I$ has two components, the points of $I$ on one branch, the points of $I^c$ on the other. Since there is an obvious relation $D_I=D_{I^c}$, we make the convention $|I \cap \{n-2,n-1,n\}| \leq 1$ to avoid this duplication. Under this assumption, it is easy to see that the inverse image of $D_I$ via the morphism $F$
is the exceptional divisor $D_I$ when $|I| \geq 3$. The case $|I|=2$ must be treated separately. The formulas are given in Remark \ref{rem}.
\end{rem}
 
\section{The Chow ring of $\overline{M}_{0,n}$}

We first review some general facts from [FM] about the intersection ring of the blow-up $\widetilde{Y}$ of a smooth variety $Y$ along a smooth irreducible subvariety $Z$. 
When the restriction map from $A^*(Y)$ to $A^*(Z)$ is surjective, S. Keel has shown in [K] that the computations become simpler. 
We denote the kernel of the restriction map by $J_{Z/Y}$ so that $$A^*(Z)=\frac{A^*(Y)}{J_{Z/Y}}.$$ Define a Chern polynomial for $Z \subset Y$, denoted by 
$P_{Z/Y}(t)$, to be a polynomial $$P_{Z/Y}(t)=t^d+a_1t^{d-1}+\etc+a_{d-1}t+a_d \in A^*(Y)[t],$$ 
where $d$ is the codimension of $Z$ in $Y$ 
and $a_i \in A^i(Y)$ is a class whose restriction in $A^i(Z)$ is $c_i(N_{Z/Y})$, where $N_{Z/Y}$ is the normal bundle of $Z$ in $Y$. We also require that $a_d=[Z]$, 
while the other classes $a_i$, for $0<i<d$, are determined only modulo $J_{Z/Y}$.

We identify $A^*(Y)$ as a subring of $A^*(\widetilde{Y})$ by means of the map $\pi^*:A^*(Y) \rightarrow A^*(\widetilde{Y})$, 
where $\pi:\widetilde{Y}\rightarrow Y$ is the birational morphism. Let $E \subset \widetilde{Y}$ be the exceptional divisor. The formula of Keel is as follows:
The Chow ring $A^*(\widetilde{Y})$ is given by 
\begin{equation}  \tag{1}\label{keel}A^*(\widetilde{Y})=\frac{A^*(Y)[E]}{(J_{Z/Y} \cdot E,P_{Z/Y}(-E))}.\end{equation}

The following facts are used to find the relation between the ideals $J_{V/Y}$ and $J_{\widetilde{V}/\widetilde{Y}}$ for a subvariety $V$ of $Y$:

\begin{itemize}
\item Suppose that $V$ is a nonsingular subvariety of $Y$ that intersects $Z$ transversally in an irreducible variety $V \cap Z$, and that the restriction $A^*(V) \rightarrow A^*(V \cap Z)$ is also surjective. 
Let $\widetilde{V}=Bl_ZV.$ Then $A^*(\widetilde{Y}) \rightarrow A^*(\widetilde{V})$ is surjective, with kernel $J_{V/Y}$ if $V \cap Z$ is not empty, and kernel $(J_{V/Y},E)$ if $V \cap Z$ is empty.

\item Suppose that $Z$ is the transversal intersection of nonsingular subvarieties $V$ and $W$ of $Y$, and that the restrictions $A^*(Y) \rightarrow A^*(V)$ and $A^*(Y) \rightarrow A^*(W)$ are also surjective. 
Let $\widetilde{V}=Bl_Z V.$ Then
\begin{enumerate}
\item $A^*(\widetilde{Y}) \rightarrow A^*(\widetilde{V})$ 
is surjective, with kernel $(J_{V/Y},P_{W/Y}(-E))$;
\item $A^*(\widetilde{Y})\rightarrow A^*(E \cap \widetilde{V})$ 
is surjective, with kernel $(J_{Z/Y},P_{W/Y}(-E))$.
\end{enumerate}
\end{itemize}

It follows that the restriction map to a blow-up center in the construction of $\overline{M}_{0,n}$ is a surjection: The inclusion of the subvariety $i_I:X_I \rightarrow X$ is a section
of a projection from $X$ onto $X_I$. This gives the injectivity of push-forward and the surjectivity of pull-back via the inclusion $i_I$. 

In our situation it is easy to determine a Chern polynomial for the blow-up centers. Since the involved blow-up centers are transversal intersection
of divisors, it will be enough to know the following two simple facts:

\begin{itemize}
\item
A Chern polynomial for a divisor class $D$ is $t+D$. 

\item
If the subvariety $Z$ is a transversal intersection of the subvarieties $V$ and $W$,  then $Z$ has a Chern polynomial $$P_{Z/Y}(t)=P_{V/Y}(t) \cdot P_{W/Y}(t).$$
\end{itemize}

Other useful formulas give the relation between a Chern polynomial $P_{V/Y}(t)$ of a subvariety $V$ and its proper transform $\widetilde{V}$:

\begin{itemize}
\item  If $V$ meets $Z$ transversally, then $P_{V/Y}(t)$ is a Chern polynomial for $\widetilde{V}$ in $\widetilde{Y}$.

\item  If $V$ contains $Z$, then $P_{V/Y}(t-E)$ is a Chern polynomial for $\widetilde{V} \subset \widetilde{Y}$.
\end{itemize}

For the proofs, see Section 5 in [FM]. 
In the study of the intersection pairing of a blow-up, there is a simple relation between certain intersection numbers in $A^*(\widetilde{Y})$ and numbers in $A^*(Y)$:

\begin{lemma}\label{num1}
Suppose that $Z$ is the transversal intersection $D_1\cap \etc \cap D_r$ of divisor classes $D_1,\etc,D_r$ on $Y$ and let $f \in A^*(Y)$ be an element of degree $d=\dim(Z)$. 
The following relation holds in $A^*(Bl_Z Y):$ $$f \cdot E^r=(-1)^{r-1}f \cdot Z.$$\end{lemma}
\begin{proof}
This is Lemma 5.3 in [T].
\end{proof}

The following vanishing result will be useful:

\begin{prop}
Let $I,J \subset \{1,\etc,n\}$ be subsets satisfying condition $(*)$. The product $D_I \cdot D_J$ is zero unless $I \subseteq J$ or $J \subseteq I$ or $I \cap J=\emptyset$.
\end{prop}
\begin{proof}
The proof is similar to that of Proposition 3.7 in [T].
\end{proof}

We are now able to study the intersection ring of the moduli space $\overline{M}_{0,n}$ more closely. As we saw above, the moduli space is constructed as a result of a sequence of blow-ups of the algebraic variety $X$. 
The Chow ring of $X$ is generated by the divisor classes $a_i$ for $1 \leq i \leq n-3$, where the class $a_i$ of degree one is defined to be the pull-back of the positive generator for the Chow ring of the projective line $\mathbb{P}^1$
via the projection onto the $i^{th}$ factor. The space of relations is generated by the relations $a_i^2=0$. The pull-backs of these classes to the moduli space are denoted by the same letters. 
 
\subsection{The space of relations} The discussion above shows that the Chow ring of the moduli space $\overline{M}_{0,n}$ is generated over $\mathbb{Z}$ by the classes $a_i$ and the exceptional divisors $D_I$.
To give a complete description of the intersection ring we need to understand the space of relations among the generators:

\begin{itemize}\label{relations}
\item There are the relations $a_i^2=0$ for $1 \leq i \leq n-3$.
\item 
For subsets $I,J \subset \{1,\etc,n\}$, satisfying condition $(*)$, the product $D_I \cdot D_J$ is zero unless $$(**) \qquad I \subseteq J, \qquad  \mathrm{or} \ J \subseteq I, \qquad \mathrm{or} \ I \cap J=\emptyset.$$

\item
For any subset $I \subset \{1,\etc,n\},$ fulfilling the requirement $(*)$, consider the inclusion $$i_I:X_I \rightarrow X.$$ 
The relations $$\{x \cdot D_I=0 : x \in \ker(i_I^*:A^*(X) \rightarrow A^*(X_I))\}$$ hold. The injectivity of the push-forward map on the Chow groups shows that this kernel is the same as the kernel of the operator 
$X_I \cap -: A^*(X) \rightarrow A^*(X)$ defined by the rule $x \rightarrow X_I \cap x$, which is easy to compute. 

\item As we saw above, in blowing up the variety $Y$ along a subvariety $Z \subset Y$, if the center $Z$ 
can be written as the transversal intersection of the subvarieties $V$ and $W$ of $Y$, then the class $P_{W/Y}(-E_Z)$ is in the ideal $J_{\widetilde{V}/\widetilde{Y}}$. 
This means that the product $P_{W/Y}(-E_Z) \cdot E_{\widetilde{V}}$ is zero, where $E_{\widetilde{V}}$ is the class of the exceptional divisor of the blow-up along the subvariety $\widetilde{V}$. 
We get a class of relations of this type by writing the centers of blow-ups introduced in the construction of the space $\overline{M}_{0,n}$ as  transversal intersections in different ways. 
If the subvariety $V$ can be written as a transversal intersection $V_1 \cap \etc \cap V_k$, we obtain the relation $P_{W/Y}(-E_Z) \cdot E_{V_1} \dots E_{V_k}=0$.

\item For each subvariety $Z \subset Y$ with a Chern polynomial $P_{Z/Y}(t)$, there is a relation $$P_{Z/Y}(-E_Z)=0,$$ where $E_Z$ is the class of the exceptional divisor of the blow-up of 
$Y$ along $Z.$ These give another class of relations in $A^*(\overline{M}_{0,n}).$
\end{itemize}

\subsection{Example} Let $I$ be a subset of the set $\{1, \etc,n-3\}$ containing at least 3 elements and let $i,j \in I$ be distinct elements. Notice that $$a_i-a_j \in \ker(i^*: A^*(X) \rightarrow A^*(X_I)),$$
from which the relation $(a_i-a_j)\cdot D_I=0$ follows. 

\subsection{Standard monomials} 
Any monomial $v \in A^d(\overline{M}_{0,n})$ can be written as a product $a(v) \cdot D(v)$, where $a(v)$ is a product of $a_i$'s and $D(v)$ is a product of exceptional divisors. The relations described in \ref{relations}
can be used to get a smaller set of generators for the intersection ring. To enumerate the monomials we associate a graph to the generator $v$. 
We first define an ordering on the \emph{polynomial ring} 
$$R:=\mathbb{Z}[a_i,D_I: 1 \leq i \leq n-3, \ I \subset \{1, \etc ,n\}, \ \mathrm{where} \  I \ \mathrm{satisfies} \ (*) ].$$

\begin{defn}\label{<}
Let $I,J \subset \{1,\etc,n\}$, we say that $I<J$ if
\begin{itemize}
\item $|I|<|J|$

\item or if $|I|=|J|$ and the smallest element in $I-I \cap J$ is smaller than the smallest element of $J-I\cap J$.
\end{itemize}
Put an arbitrary total order on monomials in $$\mathbb{Z}[a_i:1 \leq i \leq n-3].$$
Suppose $v_1,v_2 \in R$ are monomials. We say that $v_1 < v_2$ if we can write them as $$v_1=a(v_1) \cdot \prod_{r=1}^{r_0}D_{I_r}^{i_r} \cdot D, \qquad \ \mathrm{and} \ v_2=a(v_2) \cdot \prod_{r=1}^{r_0}D_{I_r}^{j_r} \cdot D,$$ 
where $D=\prod_{r=r_0+1}^m D_{I_r}^{i_r}$, for $I_m < \etc <I_1$, and $i_{r_0}<j_{r_0};$ 
or if $r_0=0$ and $a(v_1)<a(v_2)$. 

Furthermore, we say that  $v_1 \ll v_2,$ if for any factor $D_I$ of $v_2$ we have that $v_1 <D_I.$ Note that $v_1 \ll v_2$ implies that $v_1 <v_2$. 
\end{defn}

\begin{defn}\label{Graph}
Let $v=a(v)\cdot D_{I_1}^{i_1} \dots D_{I_m}^{i_m}$, where  $i_r \neq 0$ for  $r=1, \etc,m$ and $I_m< \etc< I_1$, be a monomial. 
The directed graph $\mathcal{G}=(V_\mathcal{G},E_\mathcal{G})$ associated to $v$ or of the collection $I_1,\etc,I_m$ is defined by the following data:
\begin{itemize}
\item A set $V_\mathcal{G}$ and a one-to-one correspondence between members of $V_{\mathcal{G}}$ and members of the set $\{1,\etc,m\}$. The elements of $V_\mathcal{G}$ are called the vertices of $\mathcal{G}$.
\item  A set $E_\mathcal{G} \subset V_\mathcal{G} \times V_\mathcal{G}$ consisting of all pairs $(r,s)$, where $I_s$ is a maximal element of the set $$\{I_i: I_i \subset I_r\}$$ 
with respect to inclusion. The elements of $E_\mathcal{G}$ are called the edges of $\mathcal{G}$.
\end{itemize}

For a vertex $i \in V_\mathcal{G}$, the closure $\overline{i} \subset V_\mathcal{G}$ is defined to be the subset
$$\{r \in V_\mathcal{G}: I_r \subseteq I_i\}$$ of $V_\mathcal{G}.$ The degree $\deg(i)$ of $i$ is defined to be the number of the elements of the set 
$$\{j \in V_\mathcal{G}: (i,j) \in E_\mathcal{G}\}.$$

A vertex $i \in V_\mathcal{G}$ is called a \emph{root} of $\mathcal{G}$ if $I_i$ is maximal with respect to inclusion of sets. 
Maximal vertices of $\mathcal{G}$ are called \emph{external} and all the other vertices will be called \emph{internal}.

In the following, we use the letters $I_1,\etc,I_m$ to denote the vertices of $\mathcal{G}$.
\end{defn}

\begin{rem}
It is easy to see that the graph $\mathcal{G}$ associated to a non-zero monomial $v$ has no loop. Therefore, we refer to $\mathcal{G}$ as the associated forest of $v$ or of the collection $I_1,\etc,I_m$.
\end{rem}

\begin{defn}\label{standard}
Let $v$ be as in Definition \ref{Graph}, $\mathcal{G}$ be the associated forest, and $J_1,\etc,J_s,$ for some $s \leq m$, be roots of $\mathcal{G}$ satisfying $\{n-2,n-1,n\} \cap J_r=\emptyset$. For each such root $J_r$ let $\alpha_r \in J_r$ 
be the smallest element. The subset $S$ of the set $\{1,\etc,n-3\}$ is defined as follows:
$$S:=\{\alpha_1,\etc,\alpha_s\} \cup (\cap_{r=1}^m I_r^c - \{n-2,n-1,n\}).$$

The monomial $v$ is said to be standard if
\begin{itemize}
\item The monomial $a(v) \in \mathbb{Z}[a_i: i \in S]/(a_i^2:i \in S)$ is non-zero.
\item For each $r$ we have that $$i_r  \leq \min(|I_r|-2, |I_r|-|\cup_{I_s \subset I_r} I_s|+\deg(I_r)-2).$$
\end{itemize}
\end{defn}

It is now easy to see by the same argument as in Proposition 5.9 in [T] that:

\begin{prop}\label{st}
The Chow group $A^d(\overline{M}_{0,n})$ of $\overline{M}_{0,n}$, for $0 \leq d \leq n-3$, is additively generated by standard monomials.
\end{prop}

\subsection{Definition of the dual element}

In this section we show that there is an involution on the intersection ring of $\overline{M}_{0,n}$ which gives a one-to-one correspondence between standard monomials in $A^d(\overline{M}_{0,n})$ and those in $A^{n-3-d}(\overline{M}_{0,n})$,
for $0 \leq d \leq n-3$. 

\begin{defn}\label{dual}
Suppose $v=a(v) \cdot D(v)$ is a standard monomial, where $a(v)$ is a product of $a_i$'s for $i\in A_v \subseteq \{1,\etc,n-3\}$ and $$D(v)=\prod_{r=1}^m D_{I_r}^{i_r},$$ 
where $i_r \neq 0$ for $r=1,\etc,m,$ and $I_m < \etc<I_1$. Let $\mathcal{G}$ be the associated forest, and $J_1,\etc,J_s$, the set $S$ and $\alpha_r \in J_r$ for $1 \leq r \leq s$ be as in Definition \ref{standard}. 
The subset $T$ of the set $S$ is defined to be $$T:=S - A_v.$$

For each $1 \leq r \leq m,$ define $j_r$ to be $$ \left\{ \begin{array}{ll}
|I_r|-|\cup_{I_s \subset I_r} I_s|+\deg(I_r)-1-i_r & \qquad \mathrm{if} \ I_r \ \mathrm{is \ an \ internal \ vertex \ of} \ \mathcal{G}; \\ \\
|I_r|-1-i_r & \qquad \mathrm{if} \ I_r \ \mathrm{is \ an \ external \ vertex \ of} \ \mathcal{G}.  \\
\end{array} \right. $$

We define $v^*=a(v^*) \cdot D(v^*),$ where  $$a(v^*)=\prod_{i \in T}a_i,  \qquad D(v^*)=\prod_{r=1}^m D_{I_r}^{j_r}.$$
\end{defn}

\begin{rem}
It is straightforward to see that the dual of a standard monomial $v$ is well-defined and furthermore, it is standard. The property $v^{**}=v$ shows that the operator $*$ is an involution.
\end{rem}

The next lemma will be useful in the proof of Proposition \ref{tri} and identity \eqref{Number}:

\begin{lemma}\label{deg}
Let $v=a(v) \cdot D(v)$ be as in Definition \ref{dual}, and $\mathcal{G}$ be the associated forest. For a vertex $i \in V_\mathcal{G}$ corresponding to the subset $I_i$ of the set $\{1,\etc,n\}$, 
the equality $$\sum_{\overline{i}}(i_r+j_r)=|I_i|-1$$ holds. Here $\overline{i}$ is the closure of $i$ in $\mathcal{G}$ as in Definition \ref{Graph}.
\end{lemma}
\begin{proof}
It is immediate from the definition of the $j_r$'s above.
\end{proof}

\subsection{The pairing $A^d(\overline{M}_{0,n}) \times A^{n-3-d}(\overline{M}_{0,n})$}
In the previous part, we defined dual elements for standard monomials. Below, we will see that with respect to the ordering of the generators of the Chow groups given in \ref{<} the resulting intersection matrix between the standard 
monomials and their duals consists of identity matrices along the main diagonal, up to a sign, and all blocks below the diagonal are zero. 
To prove the stated properties of the intersection matrix, we introduce a natural filtration on the Chow ring.

\begin{defn}\label{filter}
Let $v$ be a standard monomial as given in Definition \ref{standard}, and let $J_1,\etc,J_s$ be the roots of the associated forest. Define $p(v)$ to be the degree of the element $$a(v) \cap_{r=1}^s X_{J_r} \in A^*(X),$$ 
which is the same as the integer $$\deg{a(v)}+\sum_{r=1}^s |J_r|-s.$$ The subspace $F^p A^*(\overline{M}_{0,n})$ 
of the Chow ring is defined to be the $\mathbb{Z}$-submodule generated by standard monomials $v$ satisfying $p(v) \geq p$.
\end{defn}

\begin{prop}\label{fil}
(a) For any integer $p$, we have that $F^{p+1}A^*(\overline{M}_{0,n}) \subseteq F^pA^*(\overline{M}_{0,n})$.

(b) Let $v \in F^pA^*(\overline{M}_{0,n})$ and $w \in A^d(\overline{M}_{0,n})$ be such that $w \ll v.$ If $p+d > n-3,$ then $v \cdot w$ is zero. In particular, $F^{n-2} A^*(\overline{M}_{0,n})$ is zero.
\end{prop}
\begin{proof}
The proof is similar to the proof of Proposition 5.15 in [T].
\end{proof}

Using property (b) in the previous proposition we can show the triangular property of the intersection matrix:

\begin{prop}\label{tri}
Suppose $v_1,v_2 \in A^d(\overline{M}_{0,n})$ are standard monomials satisfying $D(v_1)<D(v_2)$. Then $v_1 \cdot v_2^*=0$.
\end{prop}
\begin{proof}
The monomial $v_1 \cdot v_2^*$ can be written as a product $v \cdot w$, for $v,w \in A^*(\overline{M}_{0,n})$ satisfying the properties given in Proposition \ref{fil}. For the argument, see Proposition 5.16 in [T].
\end{proof}

To study the square blocks on the main diagonal we relate the numbers occurring in the intersection matrix to numbers in the Chow group $A_0(X)$. We first recall the relevant identity:

Let $Y$ be a blow-up of $X$ at some step in the construction of $\overline{M}_{0,n}.$ Suppose that $$V_1 \cap \etc \cap V_k \cap W=Z$$ is a transversal intersection of cycles, 
where $W=D_1 \cap \etc \cap D_r$ is a transversal intersection of divisors $D_1, \etc ,D_r \in A^1(Y)$, 
and let $f \in A^*(Y)$ be an element of degree $d=\dim(Z)$. Denote by $E_Z$ the exceptional divisor of the blow-up 
$Bl_Z Y$ of $Y$ along $Z$ and by $E_{V_1},\etc ,E_{V_k}$ those of the blow-up $\widetilde{Y}$ of $Bl_Z Y$ along the proper transform of the subvarieties $V_1, \etc ,V_k$. 
In [T] it is proven that the following identity holds in $A^{r+d+k}(Y)$:
$$f \cdot E_Z^rE_{V_1} \dots E_{V_k}=(-1)^{r-1}f \cdot W \cdot E_{V_1} \etc E_{V_k}.$$

This identity reduces the number of exceptional divisors by one for certain monomials.  
If the codimension of the subvariety $V_i$ is $r_i$ and that of $Z$ is $r_0$, then from the proven result in Lemma \ref{num1} one gets the following identity:
\begin{equation}  \tag{2}\label{number}f \cdot E_Z^r \cdot E_{V_1}^{r_1} \dots E_{V_k}^{r_k}=(-1)^{r_0-k-1}f \cdot Z.\end{equation}

We now use this identity to compute the numbers occurring on the main diagonal of the intersection matrix.  
Let $I_1, \etc ,I_m$ be subsets of the set $\{1, \etc ,n\}$ such that each pair $I_r$ and $I_s$ satisfies property $(**)$.
Let $\mathcal{G}$ be the associated forest of the collection $I_1, \etc ,I_m$ and $J_1, \etc ,J_s$, the set $S$ and $\alpha_r \in J_r$ for $1 \leq r \leq s$ be as in Definition \ref{standard}. 
Define $$D:=\prod_{r=1}^m D_{I_r}^{i_r},$$ where 
$$i_r= \left\{ \begin{array}{ll}  |I_r|-|\cup_{I_s \subset I_r} I_s|+\deg(I_r)-1 & \qquad \mathrm{if} \ I_r \ \mathrm{is \ an \ internal \ vertex \ of} \ \mathcal{G} \\ \\
|I_r|-1 & \qquad \mathrm{if} \ I_r \ \mathrm{is \ an \ external \ vertex \ of} \ \mathcal{G}.  \\
\end{array} \right.$$

Consider an element $f \in \mathbb{Z}[a_i: i \in S]/(a_i^2:i \in S)$, of degree $|\cap_{r=1}^m I_r^c|+s-3$. Then from identity \eqref{number} it follows that 
\begin{equation}  \tag{3}\label{Number}f \cdot D= (-1)^{\varepsilon} \cdot f \cdot \prod_{i \in \{1,\etc ,n-3\}-S}a_i,\end{equation}
where $\varepsilon=|\cup_{r=1}^m I_r|+\sum_{i \in V(G)}\deg(i)$, by Lemma \ref{deg}. 
It is easy to see that for any $v \in A^d(\overline{M}_{0,n})$, the product $D:=D(v) \cdot D(v^*)$ is in the form given above. We can now prove our main result:

\begin{thm}\label{main}
The Chow ring $A^*(\overline{M}_{0,n})$ of $\overline{M}_{0,n}$ is the polynomial ring over $\mathbb{Z}$ with generators $a_i$, $1 \leq i \leq n-3$,  and $D_I$ for subsets $I \subset \{1,\etc ,n\}$ satisfying $(*)$, 
modulo the ideal generated by the relations in \ref{relations}.
The involution $*$ defines a duality between the Chow groups in complementary degrees. More precisely, for each $0 \leq d \leq n-3$, the intersection pairing 
$$A^d(\overline{M}_{0,n}) \times A^{n-3-d}(\overline{M}_{0,n}) \rightarrow \mathbb{Z}$$ is perfect. 
\end{thm}
\begin{proof}
We have observed that the intersection ring is generated by the divisor classes $a_i,D_I$ over $\mathbb{Z}$, and that the relations are those described in \ref{relations}. 
To prove the non-degeneracy of the intersection pairing we show that the intersection matrix $(v_i \cdot v_j^*)$ is invertible, where the $v_i$ varies over the set of standard monomials in degree $d$. 
According to the vanishing result proven in Proposition \ref{tri} it would be enough to see the invertibility of the square block of the intersection matrix consisting of all such $v$'s having the same $D$-part. 
For a fixed monomial $$D \in \mathbb{Z}[D_I: I \subset \{1, \dots ,n\} \ \mathrm{satisfies} \ (*)],$$ denote by $A_D$ the set of all standard monomials $v$ in $A^d(\overline{M}_{0,n})$ for which $D(v)=D$. 
Let $\mathcal{G}$ be the graph associated to the monomial $D$ and define $S$ as in Definition \ref{standard}. 

From identity \eqref{Number} we see that the intersection numbers $$v_1 \cdot v_2^* \in A^{n-3}(\overline{M}_{0,n})=\mathbb{Z}$$ and $$a(v_1) \cdot a(v_2^*) \in A^{|S|}((\mathbb{P}^1)^{|S|})=\mathbb{Z}$$ differ by $(-1)^{\varepsilon}$, 
where $\varepsilon=|\cup_{r=1}^m I_r|+\sum_{i \in V(\mathcal{G})}\deg(i)$. This means that the intersection matrices $(v_i \cdot v_j^*)$ and $(a(v_i) \cdot a(v_j^*))$ for $v_i,v_j \in A_D$ are the same up to a sign. 
The latter matrix is the identity matrix for the natural choice of basis for the Chow groups of $(\mathbb{P}^1)^S$. This proves the claim.
\end{proof}

\begin{rem}\label{rem}
Let us compare our result with that of Keel in [K]. He shows that the Chow ring of $\overline{M}_{0,n}$ is generated by all divisor classes $D_I$ for which $I \subset \{1,\etc ,n\}$ satisfies $|I \cap \{n-2,n-1,n\}| \leq 1$ and $|I| \geq 2$. 
In his description the markings 1,2,3 play the rule of $n-2,n-1,n$ in ours. 
Here we replace the divisors $D_I$ for which $|I|=2$ with the divisor classes $a_i$.
They are related by the following formulas: $$a_i=\sum_{i,n-2 \in I}D_I=\sum_{i,n-1 \in I}D_I=\sum_{i,n \in I}D_I,$$
$$a_i+a_j=\sum_{i,j \in I}D_I.$$
Keel shows that the ideal of relations is generated by the following ones: For any four distinct elements $i,j,k,l \in \{1,\etc ,n\}$: 
$$\sum_{\substack{i,j \in I\\k,l \notin I}}D_I=\sum_{\substack{i,k \in I\\j,l \notin I}}D_I=\sum_{\substack{i,l \in I\\j,k \notin I}}D_I,$$
and 
$$D_I \cdot D_J=0 \qquad {\text unless} \qquad I \subseteq J, \qquad J \subseteq I, \qquad I \cap J=\emptyset.$$
It is possible to derive the relations described in \ref{relations} directly from these relations. The same is true for the vanishing in Proposition \ref{tri} and the identity \eqref{Number}. 
This leads to another proof of Theorem \ref{main} independent of our construction of the moduli space. 
Although the proofs are elementary, they are longer than ours and we don't have geometric interpretations of the relations similar to those we give in \ref{relations}, \eqref{number} and \ref{fil} according to the blow-up picture. 
\end{rem}

\section{Examples}

\begin{itemize}
\item The moduli space $\overline{M}_{0,3}=M_{0,3}$ consists of a single point and its Chow ring is isomorphic to the ring of integers $\mathbb{Z}$. 

\item $\overline{M}_{0,4}=\mathbb{P}^1$, and its intersection ring is $\frac{\mathbb{Z}[a_1]}{(a_1^2)}$, where $a_1$ is the class of a point. Note that: $$a_1=D_{1,2}=D_{1,3}=D_{1,4}.$$

\item $\overline{M}_{0,5}$ is $\mathbb{P}^1 \times \mathbb{P}^1$ blown-up at 3 points. It is a del Pezzo surface of degree 5. The exceptional divisors of the blow-ups are
$D_{1,2,3},D_{1,2,4},D_{1,2,5}$. One has the following relations:
$$a_1=D_{1,3}+D_{1,2,3}=D_{1,4}+D_{1,2,4}=D_{1,5}+D_{1,2,5},$$
$$a_2=D_{2,3}+D_{1,2,3}=D_{2,4}+D_{1,2,4}=D_{2,5}+D_{1,2,5},$$
$$a_1+a_2=D_{1,2}+D_{1,2,3}+D_{1,2,4}+D_{1,2,5},$$
$$a_i \cdot D_{1,2,3}=a_i \cdot D_{1,2,4}=a_i \cdot D_{1,2,5}=0 \qquad i=1,2,$$
$$D_{1,2,3} \cdot D_{1,2,4}=D_{1,2,3} \cdot D_{1,2,5}=D_{1,2,4} \cdot D_{1,2,5}=0,$$
$$D_{1,2,3}^2=D_{1,2,4}^2=D_{1,2,5}^2=-a_1a_2.$$

In particular, from the relations above we see directly that the intersection ring is generated by the 5 classes $a_1,a_2,D_{1,2,3},D_{1,2,4},D_{1,2,5}$.
The involution $*$ permutes the classes $a_1,a_2$ and fixes the exceptional divisors. The resulting intersection matrix of the pairing $A^1(\overline{M}_{0,5}) \times A^1(\overline{M}_{0,5})$ is:

$$\left(
\begin{array}{cccc}
I_2 & 0\\
0 &-I_3 \\\end{array}
\right).$$

\item $\overline{M}_{0,6}$ is obtained from $\mathbb{P}^1 \times \mathbb{P}^1 \times \mathbb{P}^1$ as follows: Let
$$X_{1,2,3,4}=(0,0,0), \qquad X_{1,2,3,5}=(1,1,1), \qquad X_{1,2,3,6}=(\infty,\infty,\infty),$$
$$X_{1,2,4}=0 \times 0 \times \mathbb{P}^1, \qquad X_{1,3,4}=0 \times \mathbb{P}^1 \times 0, \qquad X_{2,3,4}=\mathbb{P}^1 \times 0 \times 0,$$
$$X_{1,2,5}=1 \times 1 \times \mathbb{P}^1, \qquad X_{1,3,5}=1 \times \mathbb{P}^1 \times 1, \qquad X_{2,3,5}=\mathbb{P}^1 \times 1 \times 1,$$
$$X_{1,2,6}=\infty \times \infty \times \mathbb{P}^1, \qquad X_{1,3,6}=\infty \times \mathbb{P}^1 \times \infty, \qquad X_{2,3,6}=\mathbb{P}^1 \times \infty \times \infty,$$
$$X_{1,2,3}=\{(x,x,x): x \in \mathbb{P}^1\}.$$ 

The blow-up of the 3 points $X_{1,2,3,i}$ for $i=4,5,6$ separates 10 lines $X_I$ for which $|I|=3$. Blowing-up the proper transform of these 10 lines gives $\overline{M}_{0,6}$. 
Notice that:
$$a_1=D_{1,4}+D_{1,2,4}+D_{1,3,4}+D_{1,2,3,4}=D_{1,5}+D_{1,2,5}+D_{1,3,5}+D_{1,2,3,5}=D_{1,6}+D_{1,2,6}+D_{1,3,6}+D_{1,2,3,6},$$
$$a_2=D_{2,4}+D_{1,2,4}+D_{2,3,4}+D_{1,2,3,4}=D_{2,5}+D_{1,2,5}+D_{2,3,5}+D_{1,2,3,5}=D_{2,6}+D_{1,2,6}+D_{2,3,6}+D_{1,2,3,6},$$
$$a_3=D_{3,4}+D_{1,3,4}+D_{2,3,4}+D_{1,2,3,4}=D_{3,5}+D_{1,3,5}+D_{2,3,5}+D_{1,2,3,5}=D_{3,6}+D_{1,3,6}+D_{2,3,6}+D_{1,2,3,6}.$$

There are several classes of relations. We write a few examples:
$$a_i \cdot D_{1,2,3,4}=a_i \cdot D_{1,2,3,5}=a_i \cdot D_{1,2,3,6}=0 \qquad i=1,2,3,$$
$$D_{1,2,3,4}^3=D_{1,2,3,5}^3=D_{1,2,3,6}^3=a_1a_2a_3,$$
$$a_2 \cdot D_{2,3,4}=a_3 \cdot D_{2,3,4}=0, \qquad (D_{1,2,3,4}+D_{2,3,4})^2+a_2a_3=0.$$

The duals are defined as follows:
$$a_1^*=a_2a_3, \qquad a_2^*=a_1a_3, \qquad a_3^*=a_1a_2,$$
$$D_{1,2,3,4}^*=D_{1,2,3,4}^2, \qquad D_{1,2,3,5}^*=D_{1,2,3,5}^2, \qquad D_{1,2,3,6}^*=D_{1,2,3,6}^2,$$
$$D_{1,2,4}^*=a_3 \cdot D_{1,2,4}, \qquad D_{1,3,4}^*=a_2 \cdot D_{1,3,4}, \qquad D_{2,3,4}^*=a_1 \cdot D_{2,3,4},$$
$$D_{1,2,5}^*=a_3 \cdot D_{1,2,5}, \qquad D_{1,3,5}^*=a_2 \cdot D_{1,3,5}, \qquad D_{2,3,5}^*=a_1 \cdot D_{2,3,5},$$
$$D_{1,2,6}^*=a_3 \cdot D_{1,2,6}, \qquad D_{1,3,6}^*=a_2 \cdot D_{1,3,6}, \qquad D_{2,3,6}^*=a_1 \cdot D_{2,3,6},$$
$$D_{1,2,3}^*=a_1 \cdot D_{1,2,3}=a_2 \cdot D_{1,2,3}=a_3 \cdot D_{1,2,3}.$$ 

The resulting intersection matrix of the pairing $A^1(\overline{M}_{0,6}) \times A^2(\overline{M}_{0,6})$ becomes 
$$\left(
\begin{array}{cccc}
I_6 & 0\\
0 &-I_{10} \\\end{array}
\right).$$

\item Let $n=20$ and consider two elements $v_1=\prod_{i=1}^4 D_{I_i}$ and $v_2=a_{11} \cdot \prod_{j=1}^3 D_{J_j}$ in $A^4(\overline{M}_{0,n})$, where the subsets $I_i$ and $J_j$ of the set $\{1, \etc , 20\}$ are defined as follows: 
$$I_1=\{1, \etc , 11\}, \qquad I_2=\{1,2,3\}, \qquad I_3=\{4,5,6\}, \qquad I_4=\{7,8,9\},$$
$$J_1=\{12,13,18\}, \qquad J_2=\{14,15,19\}, \qquad J_3=\{16,17,20\}.$$

The directed graphs $\mathcal{G}_1, \mathcal{G}_2$ associated to the monomials $v_1,v_2$ are pictured below. 
Note that $\mathcal{G}_1$ has 4 vertices and 3 edges. The vertex $I_1$ is internal with degree 3 and is the only root of the graph. All the other vertices are external.
The graph $\mathcal{G}_2$ has 3 vertices and no edge. The vertices of $\mathcal{G}_2$ are all external roots.
The duals are defined by:
$$v_1^*=a_1 \cdot \prod_{i=12}^{17}a_i \cdot D_{I_1}^3 D_{I_2}D_{I_3}D_{I_4},$$ 
$$v_2^*=\prod_{i=1}^{10}a_i \cdot D_{J_1}D_{J_2}D_{J_3}.$$ 
Notice that $D(v_1) <D(v_2)$ and it follows from Proposition \ref{tri} that $v_1 \cdot v_2^*=0$, which is easy to verify directly as well.
\end{itemize}

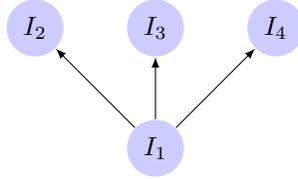
\begin{figure}[htp]
\begin{tikzpicture}
[scale=.8,auto=left,every node/.style={circle,fill=blue!20}]
\node (n1) at (3,8) {$I_1$};
\node (n2) at (1,10)  {$I_2$};
\node (n3) at (3,10)  {$I_3$};
\node (n4) at (5,10) {$I_4$};
\draw[>=latex,->] (n1) to (n2);
\draw[>=latex,->] (n1) to (n3);
\draw[>=latex,->] (n1) to (n4);
\end{tikzpicture}
\caption{The graph $\mathcal{G}_1$}
\end{figure}

\begin{figure}[htp]
\begin{tikzpicture}
[scale=.8,auto=left,every node/.style={circle,fill=blue!20}]
\node (n1) at (1,8) {$J_1$};
\node (n2) at (3,8)  {$J_2$};
\node (n3) at (5,8)  {$J_3$};
\end{tikzpicture}
\caption{The graph $\mathcal{G}_2$}
\end{figure}
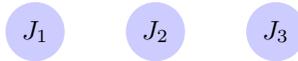

\end{document}